\documentclass[11pt,a4paper]{amsart}
\usepackage[T1]{fontenc}
\usepackage{lmodern}
\usepackage[margin=25mm]{geometry}
\usepackage{amsmath,amssymb,amsfonts,amsthm,mathtools}
\usepackage{microtype}
\usepackage{enumitem}
\usepackage[hidelinks]{hyperref}

\newtheorem{theorem}{Theorem}[section]

\newtheorem{lemma}[theorem]{Lemma}
\newtheorem{corollary}[theorem]{Corollary}

\theoremstyle{definition}

\theoremstyle{remark}
\newtheorem{remark}[theorem]{Remark}

\newcommand{\C}{\mathbb C}
\newcommand{\R}{\mathbb R}
\newcommand{\Q}{\mathbb Q}
\newcommand{\Z}{\mathbb Z}
\newcommand{\T}{\mathbb T}
\newcommand{\Pp}{\mathcal P}

\newcommand{\meas}{\operatorname{meas}}

\newcommand{\Span}{\operatorname{span}}
\newcommand{\Real}{\operatorname{Re}}

\title[Joint universality of the Riemann zeta-function]{Joint universality of the Riemann zeta-function\\ with arbitrary real shift parameters}

\begin{document}

\author{{\L}ukasz Pa\'nkowski}
\address{Faculty of Mathematics and Computer Science, Adam Mickiewicz University, Uniwersytetu Pozna\'nskiego 4, 61-614 Pozna\'{n}, Poland.}
\email{lpan@amu.edu.pl}



\begin{abstract}
We prove joint universality of the Riemann zeta-function for shifts $\zeta(s+id_j\tau)$ with arbitrary nonzero real parameters $d_j$ whose absolute values are pairwise distinct, and show that these conditions are optimal. As a consequence, we establish generalized strong recurrence for every finite family of nonzero real parameters, meaning that the corresponding shifts of the Riemann zeta-function simultaneously approximate each other. The only cases for which generalized strong recurrence remains unproved are those involving both zero and nonzero shift parameters, for which it is equivalent to the Riemann hypothesis. Our proof combines Emsalem's generalization of the six exponentials theorem with a generalization of Pecherskii's rearrangement theorem that may be of independent interest.	
\end{abstract}

\maketitle

\section{Introduction}
In 1981 Bagchi \cite{Bagchi1981} discovered an interesting connection between the Riemann hypothesis and Voronin's universality theorem \cite{Voronin1975} for the Riemann zeta-function $\zeta(s)$. Namely, he showed that the Riemann hypothesis is true if and only if, for every compact set $K\subset D:=\{s\in\C:\tfrac12<\Real s<1\}$ with connected complement and every $\varepsilon>0$, we have
\begin{equation}\label{eq:bagchi}
	\liminf_{T\to\infty}\frac1T\meas
	\left\{\tau\in[0,T]:
	\sup_{s\in K}|\zeta(s+i\tau)-\zeta(s)|
	<\varepsilon\right\}>0.
\end{equation}
In the language of topological dynamics, this is called the strong recurrence property of the Riemann zeta-function (see also \cite{Bagchi1987}). 

An interesting generalization of \eqref{eq:bagchi} was suggested by Nakamura \cite{Nakamura2009}, who considered the problem of finding all real numbers $d_1,\ldots,d_n$ for which
\begin{equation}\label{eq:nakamura-problem}
	\liminf_{T\to\infty}\frac1T\meas
	\left\{\tau\in[0,T]:
	\max_{1\le j<k\le n}\sup_{s\in K}
	|\zeta(s+id_j\tau)-\zeta(s+id_k\tau)|
	<\varepsilon\right\}>0
\end{equation}
for every such compact set $K$ and every $\varepsilon>0$. This property is usually called generalized strong recurrence or self-approximation. A natural approach to proving it is to show that the shifts $\zeta(s+id_j\tau)$ simultaneously approximate arbitrarily given non-vanishing analytic functions. Nakamura proved such a joint universality phenomenon when  $d_1,\ldots,d_n$ are algebraic real numbers linearly independent over $\Q$.

For $n=2$, the problem of joint universality for $\zeta(s+id_1\tau),\zeta(s+id_2\tau)$ was completely solved in \cite{Pankowski2009,Pankowski2016}. In \cite{Pankowski2009}, the case of irrational $d_2/d_1$ was treated using the six exponentials theorem, which shows that possible linear relations among the numbers $d_j\log p$, $j=1,2$, $p\in\mathbb{P}$, involve at most finitely many exceptional primes. The remaining case of nonzero rational $d_2/d_1\ne\pm1$ was established in \cite{Pankowski2016} by modifying the classical denseness argument based on Pecherskii's rearrangement theorem. Consequently, for two nonzero real numbers $d_1,d_2$ with $d_1\ne\pm d_2$, the shifts $\zeta(s+id_1\tau)$ and $\zeta(s+id_2\tau)$ are jointly universal. This implies that the only case for which \eqref{eq:nakamura-problem} with $n=2$ remains unresolved is $d_1=1,d_2=0$, which is equivalent to the Riemann hypothesis (see \cite[Theorem 1.4]{Pankowski2016}).

The situation for $n\ge3$ is much more difficult. In \cite{Pankowski2022}, the author proved joint universality of $\zeta(s+id_1\tau),\ldots,\zeta(s+id_n\tau)$ when $K$ is a singleton and $d_1,\ldots,d_n$ are distinct positive rational numbers, namely for any fixed $s=\sigma+it$ with $1/2<\sigma\le1$ and $t>0$, the values $\log\zeta(s+id_j\tau)$ can simultaneously approximate arbitrary complex numbers. Very recently, Mishou \cite{MishouIntegerManuscript} obtained joint universality for compact sets and positive integer parameters $d_1<\cdots<d_n$ satisfying certain additional separation conditions. Moreover, at a conference in Vilnius in September 2026, he announced that the separation conditions can be removed, so that joint universality holds for
arbitrary distinct positive integers $d_1,\ldots,d_n$. Nevertheless, the problem of characterizing the set of possible $d_1,\ldots,d_n$ satisfying joint universality of $\zeta(s+id_1\tau),\ldots,\zeta(s+id_n\tau)$, or at least \eqref{eq:nakamura-problem}, remained open. 

The main difficulty is related to the linear dependence over $\Q$ of the numbers $d_j\log p$, where $p$ runs over the primes. In the case considered by Nakamura, these numbers are linearly independent, and the classical Kronecker--Weyl theorem can be applied directly. For general real $d_j$, relations can occur both among the numbers $d_1,\ldots,d_n$ and between $d_j\log p$ corresponding to different primes. The former can be taken into account by choosing a $\Z$-basis $\beta_1,\ldots,\beta_m$ of the additive group generated by $d_1,\ldots,d_n$. To deal with the latter, we use a result due to Emsalem concerning linear relations among logarithms of algebraic numbers, which generalizes the six exponentials theorem (see \cite{WaldschmidtLectures1992,Waldschmidt2000}). As in the case of two shifts considered in \cite{Pankowski2009}, this allows us to isolate a finite set $S$ of exceptional primes such that the numbers $\beta_\ell\log p$, $1\le\ell\le m$, $p\notin S$, are linearly independent over $\Q$.

Another crucial ingredient of our proof is a generalization
of Pecherskii's rearrangement theorem.
Let $H$ be a complex Hilbert space.
In Voronin's original approach \cite{Voronin1975},
one first fixes a sequence $(v_p)$ in $H$ and
requires that, for every $0\ne g\in H$,
\[
\sum_{\substack{p}}
\max(0,\Real\langle v_p,g\rangle)=+\infty,
\qquad
\sum_{\substack{p}}
\min(0,\Real\langle v_p,g\rangle)=-\infty.
\]
Together with $\sum_p\|v_p\|^2<\infty$,
these conditions allow one to approximate
any element of $H$ by rearranging the series
$\sum_p v_p$.

On the other hand, in Bagchi's approach (see \cite{Bagchi1981} or \cite[Chapter~6, \S\,6.1, Theorem~1.16]{Laurincikas1996}), one considers vectors $z_pv_p$, where $|z_p|=1$ and $v_p\in H$ is fixed for each $p$.
The corresponding condition is
\[
\sum_p\sup_{z\in\T}
\Real\langle zv_p,g\rangle
=
\sum_p|\langle v_p,g\rangle|
=\infty
\qquad(0\ne g\in H).
\]
For example, in the proof of joint universality for Dirichlet $L$-functions, one takes $v_p=(\chi_1(p)p^{-s},\ldots,\chi_n(p)p^{-s})$, where $\chi_1,\ldots,\chi_n$ are non-equivalent Dirichlet characters.

We generalize this criterion by replacing the sets $\{zv_p:z\in\T\}$ with arbitrary nonempty bounded sets $A_p\subset H$.
More precisely, assuming
\[
\sum_p\sup_{a\in A_p}\|a\|^2<\infty,
\]
we prove that the condition
\[
\sum_p\sup_{a\in A_p}
\Real\langle a,g\rangle=\infty
\qquad(0\ne g\in H)
\]
implies that the finite sums $\sum_{p\in M}a_p$, where $a_p\in A_p$, are dense in $H$.

Our main result shows that joint universality of $\zeta(s+id_1\tau),\ldots,\zeta(s+id_n\tau)$ holds under the optimal conditions $d_j\ne0$ and $|d_j|\ne|d_k|$ for $j\ne k$. In particular, no linear independence assumptions
on $d_1,\ldots,d_n$ are needed.

\begin{theorem}\label{thm:main} Let $n\ge1$ and let $d_1,\ldots,d_n$ be nonzero real numbers satisfying 
	\begin{equation}\label{eq:distinct}
		|d_j|\ne|d_k|\qquad(1\le j<k\le n).
	\end{equation}
	For each $1\le j\le n$, let $K_j\subset D$ be compact with connected complement, and let $f_j$ be continuous and non-vanishing on $K_j$ and holomorphic in its interior. Then for every $\varepsilon>0$,
	\begin{equation}\label{eq:main}
		\liminf_{T\to\infty}\frac1T\meas
		\left\{\tau\in[T,2T]:
		\max_{1\le j\le n}\sup_{s\in K_j}
		|\zeta(s+id_j\tau)-f_j(s)|<\varepsilon
		\right\}>0.
	\end{equation}
\end{theorem}

The conditions on $d_1,\ldots,d_n$ in Theorem~\ref{thm:main} are optimal for joint universality. Indeed, if $d_j=0$, the corresponding function does not depend on $\tau$. If $d_j=d_k$, the two corresponding shifts are identical. Finally, if $d_j=-d_k$, then for every real $\sigma\in(1/2,1)$,
\[
\zeta(\sigma+id_j\tau)
=\overline{\zeta(\sigma+id_k\tau)}
\qquad(\tau\in\R),
\]
so the values of $\zeta(s+id_j\tau)$ and $\zeta(s+id_k\tau)$ cannot approximate two independent functions. 

As an immediate consequence of Theorem~\ref{thm:main}, we obtain generalized strong recurrence for every nonzero real parameters $d_1,\ldots,d_n$.

\begin{corollary}\label{cor:recurrence}
	Let $n\ge2$ and let $d_1,\ldots,d_n$ be
	arbitrary nonzero real numbers. Then for every
	compact set $K\subset D$ with connected
	complement and every $\varepsilon>0$, we have
	\begin{equation}\label{eq:allrecurrence}
		\liminf_{T\to\infty}\frac1T\meas
		\left\{\tau\in[T,2T]:
		\max_{1\le j<k\le n}\sup_{s\in K}
		|\zeta(s+id_j\tau)-\zeta(s+id_k\tau)|
		<\varepsilon
		\right\}>0.
	\end{equation}
\end{corollary}
\begin{proof}
	Let $a_1,\ldots,a_m$ be the distinct numbers
	among $|d_1|,\ldots,|d_n|$, and put
	$K^*=\{\overline s:s\in K\}$.
	Choose a closed rectangle $K_0\subset D$
	containing $K\cup K^*$.
	By Theorem~\ref{thm:main}, applied to
	the positive parameters $a_1,\ldots,a_m$
	with all approximated functions equal to $1$ on $K_0$,
	there is a set of $\tau$ of positive lower
	density such that
	\[
	\max_{1\le h\le m}\sup_{s\in K_0}
	|\zeta(s+ia_h\tau)-1|
	<\frac{\varepsilon}{2}.
	\]
	Since
	$\zeta(\overline s)=\overline{\zeta(s)}$,
	the same estimate holds on $K$ for every
	shift $\zeta(s+id_j\tau)$, including those
	with $d_j<0$. This implies \eqref{eq:allrecurrence}.
\end{proof}

Corollary \ref{cor:recurrence} states that Nakamura's generalized strong recurrence problem is settled for any finite family of nonzero real parameters. Let us emphasize that if at least one parameter $d_j$ is zero and at least one is nonzero, \eqref{eq:nakamura-problem} is equivalent to the Riemann hypothesis by Bagchi's theorem and Theorem~\ref{thm:main}.

\section{Linear independence and Emsalem's theorem}

As can be seen from the proof of the main result of \cite{Nakamura2009}, the proof of Theorem \ref{thm:main} is quite standard if the set of numbers $\{d_j\log p: j=1,2,\ldots,n, p\in\mathbb{P}\}$ is linearly independent over $\Q$. In such a case one can use the Kronecker-Weyl theorem to show that the quantities
\[
-\frac{d_j\tau\log p}{2\pi}\pmod 1,
\qquad 1\le j\le n,\quad p\in Q,
\]
can simultaneously approximate arbitrary given real numbers  $\pmod 1$, for any finite set of primes $Q$.
This is a crucial step in many known proofs of universality and, simultaneously, the main obstacle in proving Theorem~\ref{thm:main} for $d_1,\ldots,d_n$ when the numbers $d_j\log p$ are linearly dependent over $\Q$. 

The case when $\sum_{j=1}^n c_jd_j=0$ for some numbers $c_j\in\Q$, not all zero, which is equivalent to  $\sum_{j=1}^n c_jd_j\log p=0$ for every prime $p$, we shall treat by choosing
a $\mathbb Z$-basis $\beta_1,\ldots,\beta_m$ of the additive group generated by $d_1,\ldots,d_n$, as in \eqref{eq:integerbasis} below. Then, for each fixed prime $p$, the numbers
$\beta_1\log p,\ldots,\beta_m\log p$ are linearly independent over $\mathbb Q$, and every $d_j\log p$ is a linear combination of them with integer coefficients.

However, linear relations involving different primes may still occur. For example, if $d_1=\log 2$ and $d_2=\log 3$, then $d_1\log 3-d_2\log 2=0$. Such relations for $n=2$ were investigated in \cite{Pankowski2009} using the six exponentials theorem.

Our first main observation is that, by applying a generalization of the six exponentials theorem due to Emsalem (see Theorem \ref{thm:emsalem}), one can show that only finitely many primes can occur in nontrivial linear relations among the numbers
\[
\beta_\ell\log p,
\qquad 1\le\ell\le m,\quad p\in\mathbb P.
\]
More precisely, we prove that there exists a finite set of primes $S$ such that  if 
\[
\sum_p\sum_{\ell=1}^m
c_{\ell,p}\beta_\ell\log p=0,
\qquad c_{\ell,p}\in\Q,
\]
and only finitely many $c_{\ell,p}$'s are non-zero, then $c_{\ell,p}=0$ for $p\notin S$, $1\le\ell\le m$. Thus, the numbers $\beta_\ell\log p$ with $p\notin S$ are linearly independent over $\mathbb Q$ and are also independent of all the numbers $\beta_\ell\log q$ with $q\in S$.

The set $\{\beta_j\log q: q\in S,\ j=1,2,\ldots,m\}$ may still contain linear dependencies. We therefore choose a $\mathbb Z$-basis $\lambda_1,\ldots,\lambda_r$ of the additive group it generates. Consequently, the numbers
\[
\lambda_1,\ldots,\lambda_r,
\qquad
\beta_\ell\log p
\quad(p\notin S,\ 1\le\ell\le m)
\]
are linearly independent over $\mathbb Q$, which enables us to apply the Kronecker-Weyl theorem.

\subsection{Exceptional set and Emsalem's theorem}

Let
\[
 L_0=\Span_{\Q}\{\log p:p\text{ prime}\}
 \subset \{\log a:a\in\overline\Q_{>0}\}:=\mathcal{L}.
\]

\begin{theorem}	\label{thm:emsalem}
	Let $m\ge1$ and let $V$ be a complex
	hyperplane in $\C^m$ such that
	$V\cap\Q^m=\{0\}$. Then
	\[
	\dim_{\Q}(V\cap\mathcal L^m)\le m(m-1).
	\]
	In particular, $V\cap L_0^m$ is
	of finite dimension over $\Q$.
\end{theorem}

This follows from Emsalem's theorem in the form stated in \cite[Chapter~13, p.~13-2]{WaldschmidtLectures1992} (see also
\cite[Corollary~11.10]{Waldschmidt2000}). The cited result concerns the larger $\Q$-vector space $\{z\in\C:e^z\in\overline{\Q}^{\times}\}$ of complex logarithms of nonzero algebraic numbers.

Before we apply Theorem \ref{thm:emsalem} to our case, let us introduce some notations. Assume that $d_1,\ldots,d_n\in\R\setminus\{0\}$ and put
\[
 G=\sum_{j=1}^n\Z d_j\subset\R.
\]
As a finitely generated torsion-free abelian group, $G$ has a
$\Z$-basis $\beta_1,\ldots,\beta_m$ for some $1\leq m\leq n$. In particular, $\beta_1,\ldots,\beta_m$ are  linearly independent over $\Q$ and
\begin{equation}\label{eq:integerbasis}
 d_j=\sum_{\ell=1}^m a_{j\ell}\beta_\ell,
 \qquad \text{for suitable } \ a_{j\ell}\in\Z.
\end{equation}
Define
\[
 V_{\boldsymbol\beta}
 =\left\{z\in\C^m:\sum_{\ell=1}^m\beta_\ell z_\ell=0\right\}.
\]
It is a complex hyperplane, with
$V_{\boldsymbol\beta}\cap\Q^m=\{0\}$.

\begin{theorem}\label{thm:finiteS}
	There exists a finite set $S$ of prime numbers such that for any
	finitely supported family $b_p\in G$,
	\begin{equation}\label{eq:NRmodS}
	 \sum_p b_p\log p=0
	 \quad\Longrightarrow\quad
	 b_p=0\quad(p\notin S).
	\end{equation}
\end{theorem}
\begin{proof}
Applying Theorem~\ref{thm:emsalem} to
$V_{\boldsymbol\beta}$ gives that the vector space
\[
V_{\boldsymbol\beta}\cap L_0^m
\subseteq
V_{\boldsymbol\beta}\cap\mathcal L^m
\]
has finite dimension over $\Q$, since $V_{\boldsymbol\beta}\cap\mathbb Q^m=\{0\}$ and 
\[
\dim_{\mathbb Q}(V_{\boldsymbol\beta}\cap L_0^m)
\le
\dim_{\mathbb Q}
\bigl(V_{\boldsymbol\beta}\cap\mathcal L^m\bigr)
\le m(m-1).
\]

Choose a finite basis
$v^{(1)},\ldots,v^{(r)}$ of $V_{\boldsymbol\beta}\cap L_0^m$. Each coordinate of every $v^{(a)}$ is a finite linear combination over $\Q$ of logarithms of primes.
Since the numbers $\log p$, with $p$ running over all primes, are linearly independent over $\mathbb Q$,
such a representation is unique. Thus, writing
\[
v^{(a)}_j=\sum_{p\in\mathbb P}c_{j,p}^{(a)}\log p,
\qquad c_{j,p}^{(a)}\in\mathbb Q,
\]
where only finitely many coefficients are nonzero,
we define
\[
S=\bigcup_{a,j}
\{p\in\mathbb P:c_{j,p}^{(a)}\ne0\}.
\]
Clearly, $S$ is finite. Moreover, every vector of
$V_{\boldsymbol\beta}\cap L_0^m$ has all its coordinates in
$\Span_{\Q}\{\log q:q\in S\}$.

Suppose now that $b_p\in G$ and $\sum_p b_p\log p=0$.
Write $b_p=\sum_\ell c_{\ell,p}\beta_\ell$ with
$c_{\ell,p}\in\Z$. We obtain
\[
 \sum_{\ell=1}^m\beta_\ell
 \left(\sum_p c_{\ell,p}\log p\right)=0.
\]
Therefore the vector $\left(\sum_p c_{1,p}\log p,\ldots, \sum_p c_{m,p}\log p\right)\in V_{\boldsymbol\beta}\cap L_0^m$. By the construction of $S$, every coordinate has a representation
involving only $\log q$ with $q\in S$. Since the logarithms of distinct
primes are linearly independent over $\Q$, we have $c_{\ell,p}=0$ for every $p\notin S$ and every $\ell$.
This proves \eqref{eq:NRmodS}.
\end{proof}

\subsection{Parametrization of the factors \(p^{-id_j\tau}\)}
\label{sec:kronecker}

Theorem~\ref{thm:finiteS} allows us to choose real
numbers that are linearly independent over $\Q$ and hence to apply
the classical Kronecker--Weyl theorem. However, the proof of
Theorem~\ref{thm:main} involves the factors
\[
p^{-id_j\tau},\qquad 1\le j\le n,\quad p\in Q,
\]
where $Q$ is a finite set of primes, and the numbers $d_j\log p$
need not be linearly independent over $\Q$. We therefore describe
how to parametrize these factors using the independent numbers
obtained from Theorem~\ref{thm:finiteS}.

Recall from \eqref{eq:integerbasis} that
\[
G=\bigoplus_{\ell=1}^m\Z\beta_\ell,
\qquad
d_j=\sum_{\ell=1}^m a_{j\ell}\beta_\ell,
\quad a_{j\ell}\in\Z.
\]
The group $\sum_{q\in S}G\log q$ is finitely generated and torsion-free, and hence admits a $\Z$-basis $\lambda_1,\ldots,\lambda_r$. We allow $r=0$ when $S=\varnothing$. Although the numbers $\beta_\ell\log q$, $1\le\ell\le m$, $q\in S$, generate $\sum_{q\in S}G\log q$, they may satisfy nontrivial linear relations over $\Q$ (see Remark \ref{rem:depend}). For every $q\in S$ and
$1\le j\le n$, write
\begin{equation}\label{eq:exceptionalbasis}
	d_j\log q=\sum_{h=1}^r b_{j,q,h}\lambda_h,
	\qquad b_{j,q,h}\in\Z.
\end{equation}

Let $Q\supseteq S$ be a finite set of primes. By
\eqref{eq:NRmodS}, the numbers
\begin{equation}\label{eq:independentfrequencies}
	\lambda_1,\ldots,\lambda_r,
	\qquad \beta_\ell\log p
	\quad(p\in Q\setminus S,\ 1\le\ell\le m)
\end{equation}
are linearly independent over $\Q$. Indeed, consider a linear combination over $\Q$ among the numbers in \eqref{eq:independentfrequencies}. Multiplying by a common denominator and expressing the terms involving $\lambda_h$ as integer combinations
of $d_j\log q$ with $q\in S$ show
that, by \eqref{eq:NRmodS}, the coefficient in $G$ at each $p\in Q\setminus S$ must
vanish. Since $\beta_1,\ldots,\beta_m$ are linearly independent
over $\Q$, all the coefficients at such primes are zero. The
remaining coefficients are zero as well, because
$\lambda_1,\ldots,\lambda_r$ are linearly independent over $\Q$.

Thus the Kronecker--Weyl theorem applies to the numbers in
\eqref{eq:independentfrequencies}. To recover the factors
$p^{-id_j\tau}$ from \eqref{eq:independentfrequencies}, we introduce
the following parametrization.

Put $N_Q=r+m|Q\setminus S|$ and identify
\[
[0,1]^{N_Q}=[0,1]^r\times
\prod_{p\in Q\setminus S}[0,1]^m.
\]
For $\boldsymbol u=(u_1,\ldots,u_r)\in[0,1]^r$ and
$\boldsymbol v_p=(v_{p,1},\ldots,v_{p,m})\in[0,1]^m$,
$p\in Q\setminus S$, define
\[
\Psi_Q:[0,1]^{N_Q}\longrightarrow(\T^n)^{|Q|}
\]
by
\[
\Psi_Q\bigl(\boldsymbol u,
(\boldsymbol v_p)_{p\in Q\setminus S}\bigr)
=\bigl(\psi_S(\boldsymbol u),
(\psi_p(\boldsymbol v_p))_{p\in Q\setminus S}\bigr),
\]
where
\[
\psi_S(\boldsymbol u)
=\left(
\left(e\left(\sum_{h=1}^r b_{j,q,h}u_h\right)\right)_{j=1}^n
\right)_{q\in S}
\]
and
\[
\psi_p(\boldsymbol v_p)
=\left(e\left(\sum_{\ell=1}^m
a_{j\ell}v_{p,\ell}\right)\right)_{j=1}^n
\qquad(p\in Q\setminus S).
\]
Define $H_S=\psi_S([0,1]^r)$ and
\[
H_{\boldsymbol d}
=\left\{
\left(e\left(\sum_{\ell=1}^m a_{j\ell}v_\ell\right)\right)_{j=1}^n
:\ \boldsymbol v\in[0,1]^m
\right\}.
\]
The image of $\psi_p$ is $H_{\boldsymbol d}$ for every
$p\in Q\setminus S$. Consequently,
\begin{equation}\label{eq:torusdecomp}
	H_Q:=\Psi_Q([0,1]^{N_Q})
	=H_S\times\prod_{p\in Q\setminus S}H_{\boldsymbol d}.
\end{equation}

Now let us put
\[
\boldsymbol\alpha_Q
=-\frac1{2\pi}\left(
\lambda_1,\ldots,\lambda_r,
(\beta_\ell\log p)_{
	p\in Q\setminus S,\ 1\le\ell\le m}
\right),
\]
and, for $\tau\in\R$, define
\[
\gamma_Q(\tau)
=\left((p^{-id_j\tau})_{j=1}^n\right)_{p\in Q}.
\]
Since the coefficients $a_{j\ell}$ and $b_{j,q,h}$ are
integers, \eqref{eq:integerbasis} and \eqref{eq:exceptionalbasis}
give
\[
\gamma_Q(\tau)
=\Psi_Q(\{\tau\boldsymbol\alpha_Q\})\in H_Q,
\]
where the fractional parts are taken coordinatewise.
By \eqref{eq:independentfrequencies}, the coordinates of
$\boldsymbol\alpha_Q$ are linearly independent over $\Q$.
The Kronecker--Weyl theorem therefore implies that
$\{\tau\boldsymbol\alpha_Q\}$ is uniformly distributed
modulo~$1$. In particular,
\[
H_Q=\overline{\{\gamma_Q(\tau):\tau\in\R\}}.
\]
In the proof of Theorem~\ref{thm:main}, we shall apply
\cite[Theorems~A.8.1 and~A.8.3]{KaratsubaVoronin1992}
directly to $\{\tau\boldsymbol\alpha_Q\}$.

\begin{remark}\label{rem:depend}
The following example shows that $H_S$ does not have to be a product over primes in $S$.
	
	Take $n=m=2$ and
	\[
	d_1=\beta_1=\log 2,\qquad
	d_2=\beta_2=\log 3.
	\]
	Since $d_2\log 2=d_1\log 3$ and $d_1,d_2\ne0$, we have $\{2,3\}\subseteq S$. 	Moreover,
	\[
	2^{-id_2\tau}=3^{-id_1\tau}
	\qquad(\tau\in\R).
	\]
	Since
	\[
	H_S=\overline{\{\gamma_S(\tau):\tau\in\R\}},
	\]
	every element 	$\bigl((z_{q,1},z_{q,2})\bigr)_{q\in S}$ of $H_S$ satisfies $z_{2,2}=z_{3,1}$.
	
	On the other hand, since $d_1,d_2$ are linearly independent over $\Q$, the	Kronecker--Weyl theorem implies that, for every fixed $q\in S$, the projection of $H_S$ onto the coordinates corresponding to $q$ is the whole $\T^2$. Consequently, $H_S$ cannot be the product of these projections.
\end{remark}

\section{Denseness lemma} \label{sec:joint-euler-denseness}

For a bounded smooth Jordan domain $U\subset\overline{U}\subset D$, let $B^2(U)$ denote the Bergman space of analytic functions $g:U\to\C$ satisfying
\[
\|g\|_{B^2(U)}^2:=\iint_U|g(s)|^2\,d\sigma dt<\infty.
\]
We equip $B^2(U)$ with the inner product
\[
\langle f,g\rangle=\iint_U f(s)\overline{g(s)}\,d\sigma dt.
\]
It is well known (see \cite[Theorem A.2.7]{KaratsubaVoronin1992}) that, for a compact set $K\subset U$, we have
\begin{equation*}\label{eq:bergmaneval}
	\sup_{s\in K}|g(s)|\le C_{K,U}\|g\|_{B^2(U)}.
\end{equation*}
We will work with the Hilbert space $H=B^2(U_1)\times\cdots\times B^2(U_n)$, where $U_j\subset \overline{U_j}\subset D$ are bounded smooth Jordan domains and 
\[
\langle u,g\rangle_H=
\sum_{j=1}^n\iint_{U_j}u_j(s)\overline{g_j(s)}\,d\sigma dt.
\]

Our aim is to show that, for every $v=(v_1,\ldots,v_n)\in H$ and every $\varepsilon>0$, there exist a finite set $M\subset\mathbb P\setminus S$ and real numbers $\vartheta_p$, $p\in M$, such that
\begin{equation}\label{eq:HilbertDense}
	\left\|
v-
\left(
\sum_{p\in M}\sum_{k=1}^{\infty}
\frac{e(kd_j\vartheta_p)}{kp^{ks}}
\right)_{j=1}^n
\right\|_H<\varepsilon,
\end{equation}
where $S$ denotes a finite exceptional set from Theorem~\ref{thm:finiteS}.

Choose $1/2<\sigma_1<\sigma_2<1$ such that
$\sigma_1<\Real s<\sigma_2$ on all
$\overline{U_j}$. For $p\notin S$ and
$x\in\R$, define
\begin{equation}\label{eq:up}
	u_p(x)=
	\left(-\log(1-e(d_jx)p^{-s})\right)_{j=1}^n,
	\qquad
	u_{p,j}(x;s)=
	\sum_{k\ge1}\frac{e(kd_jx)}{kp^{ks}}.
\end{equation}

A standard approach to proving \eqref{eq:HilbertDense} uses Pecherskii's rearrangement theorem (see \cite[Theorem~A.6.1]{KaratsubaVoronin1992}). To apply this theorem directly, one would first need to choose a sequence of real numbers $(\vartheta_p)_{p\notin S}$, independently of $g\in H$, such that the fixed sequence $(u_p(\vartheta_p))_{p\notin S}$ satisfies the assumptions of Pecherskii's theorem. Note that the vectors $u_p(x)$ cannot, in general, be written as $zv_p$, where $z\in\T$ and $v_p\in H$ is independent of $x$. Thus, the version of Pecherskii's theorem with unimodular coefficients mentioned in Introduction does not apply directly here.

More precisely, since $H$ may be considered as a real Hilbert space with the inner product $\Real\langle u,g\rangle_H$ and
\begin{equation}\label{eq:square-summability}
\sum_{p\notin S}
\sup_{x\in\R}\|u_p(x)\|_H^2
\ll\sum_p p^{-2\sigma_1}<\infty,
\end{equation}
in order to use Pecherskii's theorem, it suffices to check (see \cite[p. 243]{KaratsubaVoronin1992}) that, for every
$0\ne g\in H$,
\[
\sum_{p\notin S}\max(0,\Real\langle u_p(\vartheta_p),g\rangle_H)
=+\infty\qquad\text{and}\qquad \sum_{p\notin S}\min(0,\Real\langle u_p(\vartheta_p),g\rangle_H)
=-\infty.
\]

To examine this condition, define
\begin{equation}\label{eq:Deltadef}
	\Delta_j(\omega)=
	\iint_{U_j}e^{-s\omega}\overline{g_j(s)}
	\,d\sigma dt,
	\qquad
	g=(g_1,\ldots,g_n)\in H.
\end{equation}
By \eqref{eq:up}, uniformly for $x\in\R$, we have
\begin{equation}\label{eq:innerProdBySimpler}
\Real\langle u_p(x),g\rangle_H
=
\Real\sum_{j=1}^n
e(d_jx)\Delta_j(\log p)
+O_g(p^{-2\sigma_1}).
\end{equation}
Thus, a direct application of Pecherskii's theorem would require a single sequence $(\vartheta_p)$ for which the real series 
\[
\sum_{p\notin S}
\Real\sum_{j=1}^n
e(d_j\vartheta_p)\Delta_j(\log p)
\]
has both divergent positive and negative parts for every $0\ne g\in H$.

One possible approach, related to those used in \cite{MishouIntegerManuscript,Pankowski2009, Pankowski2016}, would be to choose $\vartheta_p$ so that, for suitable primes $p$,
\[
\Real\bigl(e(d_j\vartheta_p)
\Delta_j(\log p)\bigr)>0
\qquad\text{for all $j$ with }
\Delta_j(\log p)\ne0.
\]
Such a requirement is unnecessarily strong and may impose additional restrictions on $d_1,\ldots,d_n$. Instead, we shall estimate the entire sum
\[
\Real\sum_{j=1}^n
e(d_j\vartheta_p)\Delta_j(\log p)
\]
directly, without imposing separate sign conditions on its individual terms.

Since the values $e(d_1x),\ldots,e(d_nx)$ need not vary independently as $x$ ranges over $\R$, finding a single sequence $(\vartheta_p)$, independent of $g$, satisfying the required divergence conditions for every $0\ne g\in H$ is difficult.

We therefore use a more flexible variant of Pecherskii's theorem (Lemma~\ref{lem:finitechoice} below) with $A_p=\{u_p(x):x\in\R\}$. In view of \eqref{eq:square-summability}, to obtain \eqref{eq:HilbertDense} it suffices to prove that 
\[
\sum_{p\notin S}\sup_{x\in\R}
\Real\langle u_p(x),g\rangle_H
=\infty
\qquad(0\ne g\in H).
\]
By \eqref{eq:innerProdBySimpler}, this follows if 
\[
\sum_{p\notin S}\sup_{x\in\R}
\Real\sum_{j=1}^n
e(d_jx)\Delta_j(\log p)
=\infty
\qquad(0\ne g\in H),
\]
since $\sum_p p^{-2\sigma_1}<\infty$. To this end, we first prove the following lemma.

\begin{lemma}\label{lem:geometry}
	Assume that $d_j\ne0$ for all $j$ and
	$|d_j|\ne|d_k|$ for $1\le j<k\le n$.
	There exists a~constant
	$c_{\boldsymbol d}>0$ such that,
	for every $A=(A_1,\ldots,A_n)\in\C^n$,
	\begin{equation}\label{eq:geometry}
		\sup_{x\in\R}
		\Real\sum_{j=1}^n A_je(d_jx)
		\ge c_{\boldsymbol d}\|A\|_2,
		\quad
		\inf_{x\in\R}
		\Real\sum_{j=1}^n A_je(d_jx)
		\le-c_{\boldsymbol d}\|A\|_2,
	\end{equation}
	where
	$\|A\|_2=
	\left(\sum_{j=1}^n|A_j|^2\right)^{1/2}$.
\end{lemma}
\begin{proof}
Let $h_A(x)=\Real\sum_j A_je(d_jx)$. Since the numbers $\pm d_j$ are pairwise distinct and nonzero, direct
	integration of exponential functions gives
	\begin{equation}\label{eq:bohrmeans}
		\lim_{X\to\infty}\frac1X\int_0^X h_A(x)\,dx=0,
		\qquad
		\lim_{X\to\infty}\frac1X\int_0^X h_A(x)^2\,dx
		=\frac12\|A\|_2^2.
	\end{equation}
	Indeed, since
	\[
	(\operatorname{Re}z)^2
	=\frac12\left(|z|^2+\operatorname{Re}z^2\right),
	\]
	and since $d_j-d_k\ne0$ for $j\ne k$ and
	$d_j+d_k\ne0$ for all $j,k$, we have
	\[
	\begin{aligned}
		\lim_{X\to\infty}\frac1X\int_0^X h_A(x)^2\,dx
		&=\frac12\lim_{X\to\infty}\frac1X\int_0^X
		\left|\sum_{j=1}^n A_je(d_jx)\right|^2\,dx\\
		&\quad+\frac12\operatorname{Re}
		\lim_{X\to\infty}\frac1X\int_0^X
		\left(\sum_{j=1}^n A_je(d_jx)\right)^2\,dx=\frac12\sum_{j=1}^n|A_j|^2.
	\end{aligned}
	\]
	
	Let $M_A=\sup_{x\in\R}h_A(x)$ and
	$m_A=-\inf_{x\in\R}h_A(x)$. By \eqref{eq:bohrmeans}, $m_A,M_A>0$ if $A$ is nonzero. On the other hand, by the Cauchy-Schwarz inequality, we have
	$M_A,m_A\le\sum_j|A_j|\le\sqrt n\|A\|_2$.
	Moreover, since $(h_A(x)+m_A)(M_A-h_A(x))\ge0$ for every $x\in\mathbb{R}$, we have
	$h_A(x)^2\le(M_A-m_A)h_A(x)+M_A m_A$.
	Hence, using \eqref{eq:bohrmeans} gives
	$M_A m_A\ge\|A\|_2^2/2$ and, in consequence, 
	\begin{equation*}\label{eq:posneg}
		M_A\ge\frac{\|A\|_2}{2\sqrt n},
		\qquad
		m_A\ge\frac{\|A\|_2}{2\sqrt n}.
	\end{equation*}

Hence, both inequalities
in \eqref{eq:geometry} hold with
$c_{\boldsymbol d}=1/(2\sqrt n)$
when $A\ne0$.
The case $A=0$ is trivial, and the
proof is complete.
\end{proof}

Now, we prove the following variant of Pecherskii's theorem.

\begin{lemma}\label{lem:finitechoice}
Let $H$ be a complex Hilbert space and let $A_k\subset H$ be nonempty bounded sets. Put $b_k=\sup_{a\in A_k}\|a\|$. 	Assume that
	\begin{equation}\label{eq:finitechoiceconditions}
		\sum_k b_k^2<\infty,
		\qquad
		\sum_k\sup_{a\in A_k}
		\Real\langle a,g\rangle=\infty
		\quad(0\ne g\in H).
	\end{equation}
Then, for every $N\ge0$, the set
\[
\left\{
\sum_{k\in M}a_k:
M\subset\{k:k>N\}\text{ is finite},
\ \ a_k\in A_k
\right\}
\]
is dense in $H$.
\end{lemma}

The proof uses a Hahn--Banach separation argument similar to those in \cite[Theorem~4(ii)]{Good1981} and \cite[Chapter~6, \S~6.1, Theorem~1.16]{Laurincikas1996}. To pass from sums with terms in $\operatorname{conv}(A_k)$ to sums with terms in $A_k$, we use the convergence of $\sum_k\sup_{a\in A_k}\|a\|^2$. A similar argument appears in \cite[Chapter~6, \S~6.1, Lemma~1.15]{Laurincikas1996}, where sums with coefficients in the unit disc are approximated by sums with coefficients
on the unit circle.

\begin{proof}
By replacing $A_k$ with $A_k\cup\{0\}$,
we may assume that $0\in A_k$.
Indeed, the first condition in
\eqref{eq:finitechoiceconditions} remains unchanged,
while the second is preserved since the
corresponding supremum can only increase.
Moreover, the set of finite sums in the
conclusion is unchanged.

	For each $k$, put $b_k=\sup_{u\in A_k}\|u\|$.	For integers $M>L$, let
	\[
	B_{L,M}=
	\left\{
	\sum_{k=L+1}^M m_k:
	m_k\in\operatorname{conv}(A_k)
	\right\},
	\]
	where $\operatorname{conv}(A_k)$
	denotes the real convex hull of $A_k$.
	Obviously, the sets $B_{L,M}$ are convex. Moreover, the fact that 
	$0\in A_k$ implies $B_{L,M}\subseteq B_{L,M+1}$. Hence
	\[
	B_L=\overline{\bigcup_{M>L}B_{L,M}}
	\]
	is closed and convex.
	
	We claim that $B_L=H$. Otherwise, by the real
	Hahn--Banach separation theorem and the Riesz
	representation theorem, there exist $g\in H$,
	$g\ne0$, and $\alpha\in\mathbb R$ such that
	\[
	\operatorname{Re}\langle u,g\rangle
	\le\alpha
	<\operatorname{Re}\langle u_0,g\rangle
	\qquad(u\in B_L)
	\]
	for some fixed $u_0\in H\setminus B_L$. On the other hand,
	\[
	\begin{aligned}
	\sup_{u\in B_L}\Real\langle u,g\rangle =\sup_{M>L} \sum_{k=L+1}^M \sup_{a\in A_k}\Real\langle a,g\rangle =\sum_{k>L} \sup_{a\in A_k}\Real\langle a,g\rangle =\infty,
	\end{aligned}
	\]
	which gives a contradiction. Thus $B_L=H$ for every $L$.
	
	Now fix $v\in H$, $\varepsilon>0$ and $N\ge0$.
	Choose $L\ge N$ so large that
	\[
	\sum_{k>L}b_k^2<\frac{\varepsilon^2}{4}.
	\]

Since $B_L=H$, there exist $M>L$ and
$m_k\in\operatorname{conv}(A_k)$,
$L<k\le M$, such that
\[
\left\|v-\sum_{k=L+1}^M m_k\right\|
<\frac{\varepsilon}{2}.
\]
By the definition of the convex hull,
for each $L<k\le M$ we can write
\[
m_k=\sum_{r=0}^{R_k}t_{k,r}u_{k,r},
\qquad
u_{k,r}\in A_k,\quad
t_{k,r}\ge0,\quad
\sum_{r=0}^{R_k}t_{k,r}=1,
\]
where $R_k$ is finite.

	We now choose indices $r_k$ with $L<k\le M$,
	successively, so that
	\[
	\left\|\sum_{k=L+1}^M(u_{k,r_k}-m_k)\right\|^2
	\le\sum_{k=L+1}^M b_k^2.
	\]
	Indeed, suppose that $r_{L+1},\ldots,r_{k-1}$
	have already been chosen, and put
	\[
	E_{k-1}=\sum_{j=L+1}^{k-1}(u_{j,r_j}-m_j),
	\qquad E_L=0.
	\]
	Since $\sum_{r=0}^{R_k}t_{k,r}(u_{k,r}-m_k)=0$,	we have
	\[
	\begin{aligned}
		\sum_{r=0}^{R_k}t_{k,r}	\|E_{k-1}+u_{k,r}-m_k\|^2	&=\|E_{k-1}\|^2
		+\sum_{r=0}^{R_k}t_{k,r}\|u_{k,r}\|^2
		-\|m_k\|^2\le\|E_{k-1}\|^2+b_k^2.
	\end{aligned}
	\]
	As the coefficients $t_{k,r}$ are nonnegative
	and sum to $1$, there exists an index $r_k$
	such that
	\[
	\|E_{k-1}+u_{k,r_k}-m_k\|^2
	\le\|E_{k-1}\|^2+b_k^2.
	\]
	Writing $E_k=E_{k-1}+u_{k,r_k}-m_k$
	and iterating, we obtain
	\[
	\left\|\sum_{k=L+1}^M(u_{k,r_k}-m_k)\right\|^2
	=\|E_M\|^2
	\le\sum_{k=L+1}^M b_k^2
	<\frac{\varepsilon^2}{4}.
	\]
	Consequently,
	\[
	\begin{aligned}
		\left\|v-\sum_{k=L+1}^M u_{k,r_k}\right\|
		\le
		\left\|v-\sum_{k=L+1}^M m_k\right\|
		+\left\|\sum_{k=L+1}^M(u_{k,r_k}-m_k)\right\|<\frac{\varepsilon}{2}+\frac{\varepsilon}{2}=\varepsilon.
	\end{aligned}
	\]
By omitting the indices for which $u_{k,r_k}=0$, we obtain a finite sum of elements of the original sets $A_k$, with all indices greater than $N$. This completes the proof.
\end{proof}

Thus, we are ready to formulate the following denseness lemma.

\begin{theorem}\label{thm:density}
	Assume that $d_j\ne0$ for all $j$ and $|d_j|\ne|d_k|$ whenever $1\le j<k\le n$. For every $v\in H$, $\varepsilon>0$ and $y>0$, and for any real $\vartheta_p$ at the primes $p\le y$ outside $S$, there exist a finite set $M\subset\Pp\setminus S$ containing all primes $p\le y$ outside $S$, and real numbers $\vartheta_p$ for all $p\in M$, such that
	\begin{equation}\label{eq:density}
		\left\|
		\sum_{p\in M}u_p(\vartheta_p)-v
		\right\|_H<\varepsilon.
	\end{equation}
\end{theorem}

\begin{proof}[Proof of Theorem~\ref{thm:density}]
With $\sigma_1$ fixed above, put
\[
A_p=\{u_p(x):x\in\R\}
\quad(p\notin S).
\]
Uniformly in $x\in\R$, \eqref{eq:up} gives
\[
\|u_p(x)\|_H\ll p^{-\sigma_1},
\qquad
\|u_p(x)-(e(d_jx)p^{-s})_{j=1}^n\|_H
\ll p^{-2\sigma_1}.
\]
In particular,
\[
\sum_{p\notin S}
\sup_{a\in A_p}\|a\|_H^2<\infty.
\]

For $0\ne g=(g_1,\ldots,g_n)\in H$,
by \eqref{eq:Deltadef} and
Lemma~\ref{lem:geometry}, we get
\begin{align*}
	\sup_{a\in A_p}\Real\langle a,g\rangle_H
	&\ge
	\sup_{x\in\R}
	\Real\sum_{j=1}^n
	e(d_jx)\Delta_j(\log p)
	+O_g(p^{-2\sigma_1})\\
	&\ge c_{\boldsymbol d}
	\left(
	\sum_{j=1}^n
	|\Delta_j(\log p)|^2
	\right)^{1/2}
	+O_g(p^{-2\sigma_1}).
\end{align*}
	Choose $j_0$ with $g_{j_0}\ne0$. We claim that $\sum_{p\notin S}|\Delta_{j_0}(\log p)|=\infty$. 
	By \cite[Lemma 3,4, and p. 229]{KaczorowskiKulas2007}, there are $x_\nu\to\infty$, constants $B,c>0$, and
	intervals $J_\nu\subset[x_\nu,x_\nu+1]$ such that
	\begin{equation*}\label{eq:laplacecontrol}
		|J_\nu|\ge Bx_\nu^{-2},\qquad
		|\Delta_{j_0}(t)|\ge ce^{-\sigma_2x_\nu}\quad(t\in J_\nu).
	\end{equation*}
	The prime number theorem gives $|\{p:\log p\in J_\nu\}|	\gg\frac{e^{x_\nu}}{x_\nu^{3}}$. Hence, since $\sigma_2<1$,
	\[
	\sum_{\log p\in J_\nu}|\Delta_{j_0}(\log p)|
	\gg\frac{e^{(1-\sigma_2)x_\nu}}{x_\nu^{3}}
	\longrightarrow\infty,
	\]
	which proves our claim. 	
	
Thus, we have
\[
\sum_{p\notin S}
\sup_{a\in A_p}\Real\langle a,g\rangle_H
=\infty,
\]
and we are ready to apply Lemma~\ref{lem:finitechoice}. Let us fix $y>0$ and choose $y'\ge y$ such that $S\subset\{p\in\mathbb P:p\le y'\}$. Write $P_{y'}= \{p\in\mathbb P\setminus S:p\le y'\}$. 
For each $p\in P_{y'}$, retain already fixed $\vartheta_p$ for $p\le y$, and choose $\vartheta_p\in\R$ arbitrarily for $y<p\le y'$. Put $v_{y'}=v-\sum_{p\in P_{y'}}u_p(\vartheta_p)$. Then, by Lemma~\ref{lem:finitechoice}, one can find a finite set of primes $M'\subset\{p\in \mathbb{P}: p>y'\}$ and real numbers $\vartheta_p$ for each $p\in M'$ such that
\[
\left\|\sum_{p\in M'}u_p(\vartheta_{p}) - v_{y'}\right\|_H<\varepsilon.
\]
Therefore, taking $M=P_{y'}\cup M'$ yields \eqref{eq:density}.
\end{proof}

\section{Proof of Theorem~\ref{thm:main}}
\label{sec:completion}

We now complete the proof by following the standard approach originating with Voronin. The only point requiring special attention is the finite exceptional set $S$.

\begin{proof}[Proof of Theorem~\ref{thm:main}]
	Fix compact sets $K_j$, functions $f_j$, and $\varepsilon>0$ as in the statement.  By Mergelyan's theorem, since $f_j$ has no zeros on $K_j$, one can find polynomials $P_j$ such that
	\begin{equation}\label{eq:Ptargets}
		\max_{1\le j\le n}\sup_{s\in K_j}
		|e^{P_j(s)}-f_j(s)|<\varepsilon/3.
	\end{equation}
	Choose bounded smooth Jordan domains $U_j$ with
	$K_j\subset U_j\subset \overline{U_j}\subset  D$, and fix real numbers
	$1/2<\sigma_1<\sigma_2<1$ such that
	$\sigma_1<\Real s<\sigma_2$ on every $\overline U_j$.
	Write $H=B^2(U_1)\times\cdots\times B^2(U_n)$ and choose $C_B\ge1$ so that
	\begin{equation}\label{eq:evalall}
		\max_j\sup_{s\in K_j}|g_j(s)|\le C_B\|g\|_H
		\qquad(g=(g_1,\ldots,g_n)\in H).
	\end{equation}
	
	Put $B=\max_j\sup_{s\in K_j}e^{\Real P_j(s)}$, and fix
	$0<\eta<1/8$ sufficiently small that
	\begin{equation*}\label{eq:expclose}
		B\bigl((1+\eta)e^{5\eta/4}-1\bigr)<\varepsilon/3.
	\end{equation*}
	
	Moreover, there is a constant $C_0$ depending only on $U_j$ and $\sigma_1$ such that 
	\[
	\sum_{j=1}^n\iint_{U_j}
	\sum_{p>y}\sum_{k\ge1}
	\frac{1}{k^2p^{2k\sigma}}\,d\sigma dt
	\le C_0y^{1-2\sigma_1}<\frac{\eta^2}{64C_B^2}
	\]
	for sufficiently large $y\geq 2$
	
	Let $S$ be the finite exceptional set from
	Theorem~\ref{thm:finiteS}, and define
	\[
	L_S(s)=\sum_{q\in S}\log(1-q^{-s}),
	\qquad \sigma>0.
	\]
By Theorem~\ref{thm:density} and \eqref{eq:evalall}, we obtain a finite set $M\subset\mathbb P\setminus S$ containing every prime $p\le y$
	outside $S$, and real numbers $\vartheta_p$, $p\in M$, such that
	\begin{equation}\label{eq:targetfinite}
		\max_j\sup_{s\in K_j}
		\left|\sum_{p\in M}\sum_{k\geq 1}\frac{e(kd_j \vartheta_p)}{kp^{ks}}-P_j(s)-L_S(s)\right|
		<\frac\eta8.
	\end{equation}
	Set $\vartheta_p = 0$ for $p\in S$ and $C=S\cup M$.  In particular, $C$ contains all primes $p\le y$.
	
	For each $p\in M$, put
	\[
	v_{p,\ell}\equiv\beta_\ell\vartheta_p\pmod1,
	\qquad 1\le\ell\le m.
	\]
	Since $a_{j\ell}\in\mathbb Z$ and
	$d_j=\sum_{\ell=1}^m a_{j\ell}\beta_\ell$, we have
	\[
	e(d_j\vartheta_p)=e\left(\vartheta_p
	\sum_{\ell=1}^m a_{j\ell}\beta_\ell\right)=e\left(\sum_{\ell=1}^m a_{j\ell}v_{p,\ell}\right)
	\]
	for every $j=1,\ldots,n$.
	Hence the vector
	\[
	\omega^*_p = (e(d_1\vartheta_p),\ldots,e(d_n\vartheta_p))
	\]
	belongs to $H_{\boldsymbol d}$. 	For the exceptional block choose all its basis parameters to be zero.
	Hence, by \eqref{eq:torusdecomp}, we have
	\[
	\omega^* = \left(\big((\underbrace{1,\ldots,1}_{n\ \text{times}})\big)_{q\in S},(e(d_1\vartheta_p),\ldots,e(d_n\vartheta_p))_{p\in M} \right)\in H_C.
	\]

	For $\omega=(\omega_{p,j})_{p\in C,\,1\le j\le n}
	\in H_C$, define
	\[
	F_j(s,\omega)
	=
	\sum_{p\in C}\sum_{k=1}^{\infty}
	\frac{\omega_{p,j}^{\,k}}{kp^{ks}},
	\qquad 1\le j\le n.
	\]
	Since $C$ is finite and $|\omega_{p,j}|=1$,
	these series converge absolutely and
	uniformly for $s\in\overline{U_j}$.
	Moreover, let us recall (see Section \ref{sec:kronecker}) that for the
	map $\gamma_C(\tau) = \left(\left(e(-\frac{d_j\tau\log p}{2\pi})\right)_{j=1}^n\right)_{p\in C}$, we have $\gamma_C(\tau)=\Psi_C\bigl(\{\tau\boldsymbol\alpha_C\}
\bigr)\in H_C$, where
	\[
		\boldsymbol\alpha_C
	=-\frac1{2\pi}
	\left(
	\lambda_1,\ldots,\lambda_r,
	(\beta_\ell\log p)_{
		p\in M,\ 1\le\ell\le m}
	\right),
	\]
	Hence,
	\begin{equation}\label{eq:Euler-phase-identity}
		F_j(s,\gamma_C(\tau)) = \sum_{p\in C}\sum_{k\ge1}	\frac{e(-k\frac{d_j\tau\log p}{2\pi})}{kp^{ks}}.
	\end{equation}
	
On the other hand, by the definition of $\omega^*\in H_C$, we have
\[
	\begin{aligned}
		F_j(s,\omega^*)
		&=
		\sum_{q\in S}\sum_{k=1}^{\infty}
		\frac{1}{kq^{ks}}
		+
		\sum_{p\in M}\sum_{k=1}^{\infty}
		\frac{e(kd_j\vartheta_p)}{kp^{ks}}=
		-L_S(s)+\sum_{p\in M}\sum_{k=1}^{\infty}
		\frac{e(kd_j\vartheta_p)}{kp^{ks}}.
	\end{aligned}
	\]
	Therefore, by \eqref{eq:targetfinite},
	\begin{equation*}\label{eq:Euler-at-omega-star}
		\max_{1\le j\le n}
		\sup_{s\in K_j}
		|F_j(s,\omega^*)-P_j(s)|
		<\frac{\eta}{8}.
	\end{equation*}
	
	We next verify that this estimate remains
	valid for vectors sufficiently close
	to $\omega^*$. Recall that
	$\Real s\ge\sigma_1>1/2$ on each $K_j$.
	For any $z,w\in\mathbb T$ and $k\ge1$,
	the factorization
	\[
	z^k-w^k
	=(z-w)\sum_{\ell=0}^{k-1}
	z^{k-1-\ell}w^\ell
	\]
	implies
	\[
	|z^k-w^k|\le k|z-w|.
	\]
	
	
Therefore, one can easily notice that for sufficiently small $\rho>0$ and for every \[
\omega\in
	\mathcal V:=\left\{	\omega\in H_C:\max_{\substack{p\in C\\1\le j\le n}}	|\omega_{p,j}-\omega^*_{p,j}|<\rho \right\},
	\]we have
	\[
		\max_{1\le j\le n}\sup_{s\in K_j}|F_j(s,\omega) - P_j(s)|	<\frac{\eta}{4}.
	\]
	
	Finally, using \eqref{eq:Euler-phase-identity},
	we conclude that
	\begin{equation}\label{eq:controlledfinite}
		\max_{1\le j\le n}\sup_{s\in K_j}
		|F_j(s,\gamma_C(\tau))-P_j(s)|
		<\frac{\eta}{4}
	\end{equation}
	provided $\gamma_C(\tau)\in\mathcal V$.

	Recall that $\omega^*\in H_C$, so there is $\boldsymbol t^*\in[0,1)^{N_C}$
	such that $\Psi_C(\boldsymbol t^*)=\omega^*$. More precisely, we may take
	\[
	\boldsymbol t^*
	=
	\left(
	\underbrace{0,\ldots,0}_{r},
	\bigl(\{\beta_\ell\vartheta_p\}\bigr)_{
		p\in M,\ 1\le\ell\le m}
	\right).
	\]
	
	Since $\Psi_C$ is continuous and $\mathcal V$
	is an open neighborhood of $\omega^*$
	in $H_C$, there exists $0<\delta_0<1/2$
	such that
	\[
	\Psi_C(B_{\delta_0})\subset\mathcal V,
	\]
	where
	\[
	B_{\delta_0}
	=
	\left\{
	\boldsymbol t\in[0,1]^{N_C}:
	\max_{1\le\nu\le N_C}
	\|t_\nu-t_\nu^*\|	\leq\delta_0
	\right\}.
	\]
	Here $\|\cdot\|$
	denotes the distance to the nearest integer.
	In particular, the Lebesgue measure of
	$B_{\delta_0}$ is
	\[
	\operatorname{meas}(B_{\delta_0})
	=(2\delta_0)^{N_C}>0.
	\]
	
For $T>0$, define
	\[
	A_T=
	\left\{
	\tau\in[T,2T]:
	\{\tau\boldsymbol\alpha_C\}
	\in B_{\delta_0}
	\right\}.
	\]
	Then, by \eqref{eq:independentfrequencies} and the
	Kronecker--Weyl theorem (see \cite[Theorem~A.8.1]{KaratsubaVoronin1992}), we have
	\begin{equation}\label{eq:weightmass}
		\lim_{T\to\infty}
		\frac{\operatorname{meas}(A_T)}{T}
		=(2\delta_0)^{N_C}>0.
	\end{equation}
	
	Moreover, for every $\tau\in A_T$, we have $\gamma_C(\tau) =\Psi_C(\{\tau\boldsymbol\alpha_C\})	\in\mathcal V$, so \eqref{eq:controlledfinite} holds for every $\tau\in A_T$.

Let $z\ge2$ be fixed, with
$C\subseteq Q_z:=\{p:p\le z\}$.
Recalling \eqref{eq:torusdecomp} gives
\[
\boldsymbol t
=
\bigl(\boldsymbol u,
(\boldsymbol v_p)_{p\in Q_z\setminus S}\bigr)
\in[0,1]^{N_{Q_z}},
\]
where $\boldsymbol u=(u_1,\ldots,u_r)\in[0,1]^r$
parametrizes the exceptional block $H_S$, while
$\boldsymbol v_p=(v_{p,1},\ldots,v_{p,m})
\in[0,1]^m$ parametrizes $H_{\boldsymbol d}$ with $p\notin S$.

Let us put
\[
\boldsymbol t
=
\bigl(\boldsymbol t_C,
(\boldsymbol v_p)_{p\in Q_z\setminus C}\bigr)\qquad\text{with}\quad \boldsymbol t_C
=
\bigl(\boldsymbol u,(\boldsymbol v_p)_{p\in M}\bigr)
\in[0,1]^{N_C},
\]
and 
\begin{equation*}\label{eq:Rdef}
	R_{j,z}(s;\boldsymbol t)
	=
	\sum_{p\in Q_z\setminus C}\sum_{k\ge1}
	\frac{
		e(k\sum_{\ell=1}^m a_{j\ell}v_{p,\ell})}{kp^{ks}},
\end{equation*}
where
$(a_{j1},\ldots,a_{jm})
\in\mathbb Z^m\setminus\{\boldsymbol0\}$ were defined in Section \ref{sec:kronecker}.

By the definition of $\boldsymbol\alpha_{Q_z}$,
\begin{align}\label{eq:Rjz}
	R_{j,z}(s;\{\tau\boldsymbol\alpha_{Q_z}\})
	&=
	\sum_{p\in Q_z\setminus C}\sum_{k\ge1}
	\frac{1}{kp^{k(s+id_j\tau)}}.
\end{align}

Put $\mathcal D_z=B_{\delta_0}\times \prod_{p\in Q_z\setminus C}[0,1]^m$ and
\[
\mathcal F_z(\boldsymbol t)
=\sum_{j=1}^n
\|R_{j,z}(\cdot;\boldsymbol t)\|_{B^2(U_j)}^2
=\sum_{j=1}^n\iint_{U_j}
|R_{j,z}(s;\boldsymbol t)|^2\,d\sigma dt.
\]
Since the function $\mathcal F_z$ is continuous on
$[0,1]^{N_{Q_z}}$, we get by \cite[Theorem~A.8.3]{KaratsubaVoronin1992} that
\begin{align*}\label{eq:weightedidentity}
	\lim_{T\to\infty}\frac1T
	\int_{A_T}
	\sum_{j=1}^n
	\|R_{j,z}(s;\{\tau\boldsymbol\alpha_{Q_z}\})\|_{B^2(U_j)}^2
	\,d\tau\notag &=
	\int_{\mathcal D_z}
	\mathcal F_z(\boldsymbol t)\,d\boldsymbol t\notag\\
	&=
	(2\delta_0)^{N_C}
	\sum_{j=1}^n\iint_{U_j}
	\sum_{\substack{p\le z\\p\notin C}}
	\sum_{k\ge1}
	\frac{1}{k^2p^{2k\Real s}}\,d\sigma dt.
\end{align*}
The last equality follows from
$\operatorname{meas}(B_{\delta_0})
=(2\delta_0)^{N_C}$ and the orthogonality
relation
\[
\int_{\prod_{v\in Q_z\setminus C}[0,1]^m}
e\left(k_1\sum_{\ell=1}^m a_{j\ell}v_{p_1,\ell}
-k_2\sum_{\ell=1}^m a_{j\ell}v_{p_2,\ell}\right)
\prod_{v\in Q_z\setminus C}d\boldsymbol v_v
=
\begin{cases}
	1,&p_1=p_2,\ k_1=k_2,\\
	0,&\text{otherwise},
\end{cases}
\]
which holds since $\boldsymbol a_j\ne
\boldsymbol0$.

	Since $C$ contains every prime $p\le y$,
	we have
	\[
	\begin{aligned}
		&\sum_{j=1}^n\iint_{U_j}
		\sum_{\substack{p\le z\\p\notin C}}
		\sum_{k\ge1}
		\frac{1}{k^2p^{2k\Real s}}\,d\sigma dt \ll \sum_{p>y}p^{-2\sigma_1}
		\ll y^{1-2\sigma_1}.
	\end{aligned}
	\]
	Thus, with a constant $C_0$
	independent of $M$, $z$, and
	$\delta_0$, we obtain
	\[
	\begin{aligned}
		\lim_{T\to\infty}\frac1T
		\int_{A_T}
		\sum_{j=1}^n
		\|R_{j,z}(s;\{\tau\boldsymbol\alpha_{Q_z}\})\|_{B^2(U_j)}^2
		\,d\tau \le
		(2\delta_0)^{N_C}
		C_0y^{1-2\sigma_1} <
		\frac{(2\delta_0)^{N_C}\eta^2}
		{64C_B^2},
	\end{aligned}
	\]
	where the last inequality follows
	from our earlier choice of $y$.
	
	By \eqref{eq:evalall}, if
	\[
	\max_j\sup_{s\in K_j}
	|R_{j,z}(s;\{\tau\boldsymbol\alpha_{Q_z}\})|>\eta,
	\]
	then
	\[
	\sum_{j=1}^n
	\|R_{j,z}(s;\{\tau\boldsymbol\alpha_{Q_z}\})\|_{B^2(U_j)}^2
	>\frac{\eta^2}{C_B^2}.
	\]
	Consequently, Markov's inequality gives
	\begin{equation}\label{eq:Rbad}
		\begin{aligned}
			&\limsup_{T\to\infty}\frac1T
			\operatorname{meas}
			\left\{
			\tau\in A_T:
			\max_j\sup_{s\in K_j}
			|R_{j,z}(s;\{\tau\boldsymbol\alpha_{Q_z}\})|>\eta
			\right\}\le
			\frac{(2\delta_0)^{N_C}}{64}.
		\end{aligned}
	\end{equation}
Thus, it follows from \eqref{eq:controlledfinite}, \eqref{eq:Rjz} and \eqref{eq:Rbad} that
\begin{equation*}\label{eq:byFiniteProdEta}
	\max_j\sup_{s\in K_j}
\left|
\sum_{p\le z}\sum_{k\ge1}
\frac{1}{kp^{k(s+id_j\tau)}}-P_j(s)
\right|
<\frac{5\eta}{4}
\end{equation*}
for all $\tau\in A_T$ outside an exceptional
set of upper density at most
$(2\delta_0)^{N_C}/64$.

\begin{align}\label{eq:byFiniteProd}
	&\liminf_{T\to\infty}\frac1T
	\operatorname{meas}\Bigg\{
	\tau\in[T,2T]:
	\max_j\sup_{s\in K_j}
	\left|
	\prod_{p\le z}(1-p^{-s-id_j\tau})^{-1}
	-e^{P_j(s)}
	\right|<\frac{\varepsilon}{3}
	\Bigg\}\geq\frac{63}{64}(2\delta_0)^{N_C}>0,
\end{align}
so it remains to show that the Riemann
zeta-function can be approximated by
these finite Euler products.

We apply Carlson's theorem
\cite[Theorem~A.2.10]{KaratsubaVoronin1992}
to the entire function
\[
G_z(s)=(1-2^{1-s})
\left(
\zeta(s)\prod_{p\le z}(1-p^{-s})-1
\right).
\]
For $\Real s>1$, it has the Dirichlet series representation $G_z(s)=\sum_{m>z}\frac{b_z(m)}{m^s}$ with $|b_z(m)|\le2$. Hence Carlson's theorem gives
\[
\begin{aligned}
	\lim_{T\to\infty}\frac1T
	\int_0^T|G_z(\sigma+it)|^2\,dt
	&=\sum_{m>z}\frac{|b_z(m)|^2}{m^{2\sigma}}\ll\sum_{m>z}m^{-2\sigma}\ll z^{1-2\sigma_1},
\end{aligned}
\]
uniformly for $\sigma_1\le\sigma\le\sigma_2$.

For every $d_j\in\mathbb R\setminus\{0\}$ and $s=\sigma+iu\in\overline U_j$, the substitution $t=u+d_j\tau$ gives
\[
	\begin{aligned}
		\lim_{T\to\infty}\frac1T\int_T^{2T}|G_z(s+id_j\tau)|^2\,d\tau=\lim_{T\to\infty}\frac1{|d_j|T}
		\int_{\min(u+d_jT,u+2d_jT)}^{\max(u+d_jT,u+2d_jT)}
		|G_z(\sigma+it)|^2\,dt\ll z^{1-2\sigma_1}.
	\end{aligned}
\]
For $d_j<0$, we also use
$|G_z(\sigma-it)|=|G_z(\sigma+it)|$,
since the Dirichlet coefficients
$b_z(m)$ are real. Thus, 
\begin{equation}\label{eq:zetameansquare}
	\begin{aligned}
		&\limsup_{T\to\infty}\frac1T
		\int_T^{2T}\sum_{j=1}^n
		\iint_{U_j}\left|
		\zeta(s+id_j\tau)
		\prod_{p\le z}
		(1-p^{-s-id_j\tau})-1
		\right|^2\,d\sigma dt\,d\tau\le C_1z^{1-2\sigma_1},
	\end{aligned}
\end{equation}
where $C_1$ is independent of $z$. Hence, we may choose $z$ sufficiently large such that $C\subseteq Q_z$ and
\begin{equation*}\label{eq:zchoice}
	C_1z^{1-2\sigma_1}
	<
	\frac{(2\delta_0)^{N_C}\eta^2}{64C_B^2}.
\end{equation*}
By \eqref{eq:evalall}, \eqref{eq:zetameansquare},
and Markov's inequality,
\begin{equation}\label{eq:Qbad}
	\begin{aligned}
		&\limsup_{T\to\infty}\frac1T
		\operatorname{meas}
		\left\{
		\tau\in A_T:
		\max_j\sup_{s\in K_j}
		\left|\zeta(s+id_j\tau)
		\prod_{p\le z}(1-p^{-s-id_j\tau})-1\right|>\eta
		\right\}\le\frac{(2\delta_0)^{N_C}}{64}.
	\end{aligned}
\end{equation}

Let $\mathcal G_T$ be the set of
$\tau\in A_T$ for which
\[
\max_j\sup_{s\in K_j}
\left|\sum_{p\in Q_z\setminus C}\sum_{k\ge1}
\frac{1}{kp^{k(s+id_j\tau)}}\right|
\le\eta,
\qquad
\max_j\sup_{s\in K_j}\left|\zeta(s+id_j\tau)
\prod_{p\le z}(1-p^{-s-id_j\tau})-1\right|
\le\eta.
\]
By \eqref{eq:weightmass}, \eqref{eq:Rbad} and \eqref{eq:Qbad}, we have
\begin{equation*}\label{eq:positivemass}
	\liminf_{T\to\infty}
	\frac{\operatorname{meas}(\mathcal G_T)}{T}
	\ge\frac{31}{32}(2\delta_0)^{N_C}>0.
\end{equation*}

For $\tau\in\mathcal G_T$, by \eqref{eq:byFiniteProdEta}, we have
\[
\begin{aligned}
	&\max_j\sup_{s\in K_j}
	|\zeta(s+id_j\tau)-e^{P_j(s)}|\le
	\max_j\sup_{s\in K_j}e^{\Real P_j(s)}
	\left((1+\eta)e^{5\eta/4}-1\right)
	<\frac{\varepsilon}{3},
\end{aligned}
\]
by our choice of $\eta$.
Together with \eqref{eq:Ptargets}, this gives
\[
\max_j\sup_{s\in K_j}
|\zeta(s+id_j\tau)-f_j(s)|
<\varepsilon
\qquad(\tau\in\mathcal G_T).
\]
This completes the proof.
\end{proof}

\end{document}